\documentclass[11pt]{amsart}
\usepackage{amsmath,amssymb,amscd,epic,eepic}
\newcommand{\A}{{\alpha}}
\newcommand{\D}{{\delta}}

\newcommand{\B}{{\beta}}

\newcommand{\R}{{\mathbb R}}    

\newcommand{{\C}}{{\mathbb C}}    

\newcommand{\fg}{{\mathfrak g}}

\newcommand{\T}{{\mathfrak T}}

\def\LA{{\lambda}}

\def\g{{\gamma}}
\def\G{{\Gamma}}
\def\RA{{\rightarrow}}
\def\bg{Bruhat graph}
\def\tfp{$T$-fixed point}
\def\pt{Peterson translate}

\def\tx{\Theta_x(X)}
\def\TX{{\mathbb T}_x(X)}
\def\cm{Cohen-Macaulay}
\def\iff{if and only if }
\def\SV{Schubert variety}
\def\SVS{Schubert varieties}
\def\st{such that }
\def\ms{\medskip}

\def\H{{\mathcal H}}
 
  \DeclareMathOperator{\IN}{in}
  \DeclareMathOperator{\gr}{gr}

\def\Black{}

\renewcommand{\section}[1]{\refstepcounter{section}\par\bigskip
\noindent\begin{center}{\normalsize\bf \thesection
.\hspace{3mm}{#1}}\end{center}\medskip \nopagebreak}

\def\theckbibliography#1{\par\bigskip
\begin{center}
{\normalsize \bf References}
\end{center}
\par
\noindent\list
  {[\arabic{enumi}]}{\settowidth\labelwidth{[#1]}\leftmargin\labelwidth
  \advance\leftmargin\labelsep
  \usecounter{enumi}}

  \sloppy\clubpenalty4000\widowpenalty4000
  \sfcode`\.=1000\relax}

\newtheorem{theorem}{Theorem}[section]

\newtheorem{corollary}[theorem]{Corollary}

\newtheorem{lemma}[theorem]{Lemma}

\newtheorem{proposition}[theorem]{Proposition}

\newenvironment{example}{\refstepcounter{theorem}\medskip\noindent{\bf
Example \thetheorem}\hspace{1mm}}{\medskip}

\newenvironment{definition}{\refstepcounter{theorem}\medskip\noindent{\bf

Definition \thetheorem}\hspace{1mm}}{\medskip}

\newenvironment{remark}{\refstepcounter{theorem}\medskip\noindent{\bf
Remark \thetheorem}\hspace{1mm}}{\medskip}

\begin{document}
\begin{center}
\bigskip

{\Large\bf Singularities of Schubert Varieties, Tangent Cones and Bruhat 
Graphs} \\
(12-12-2003)

\bigskip
\bigskip
\bigskip{\sc James B.\ Carrell\footnote{The first author was partially
supported by the Natural Sciences and Engineering Research Council
of Canada}\bigskip
\\Jochen Kuttler\footnote{The second author was partially supported by
the SNF (Schweizerischer Nationalfonds)}}
\end{center}

\bigskip
\vspace{0.74in}
\begin{center}{\bf Abstract} 
\end{center}
{\tiny Let $G$ be a semi-simple algebraic group  over ${\mathbb C}$, $B$ a
Borel subgroup of $G$, $T$ a maximal torus in $B$ and $P$ a parabolic in $G$
containing $B$. In a previous work \cite{ck}, the
authors classified the singular \tfp s of an irreducible $T$-stable subvariety $X$ of 
the generalized flag variety $G/P$. It
turns out that under the restriction that $G$ doesn't contain any $G_2$-factors,
the key geometric invariant determining the singular \tfp s of $X$
is the linear span $\tx$ of the reduced tangent
cone to $X$ at a \tfp\ $x$. The goal of this paper is to describe this 
invariant at the maximal singular \tfp s when $X$ a \SV\ in $G/B$
and $G$ doesn't contain any $G_2$-factors. 
We first decribe $\tx$ solely in terms of \pt s,
which were the main tool in \cite{ck}. 
Then, taking a further look at the \pt s (with the $G_2$-restriction), 
we are able to describe $\tx$ in terms of its
isotropy submodule  and the \bg \ of $X$ at $x$.
This refinement gives a purely 
root theoretic description, which should be useful for computations.
Finally, still with the $G_2$-restriction, these considerations lead us to 
a non-recursive algorithm for $X$'s singular locus solely involving only
the root system of $(G,T)$ and the \bg\ of $X$. }

\section{Introduction}\label{intro}
Let $G$ be a semi-simple
algebraic group over an arbitrary algebraically
closed field $k$,
and suppose $T\subset B\subset P$ are respectively  a 
maximal torus, a Borel subgroup and an arbitrary standard parabolic in $G$.
Each $G/P$, including $G/B$,
is a  projective $G$-variety with only finitely many $B$-orbits. Every $B$-orbit
contains a unique \tfp\ $x\in (G/P)^T$, and these cells 
define an affine paving of $G/P$. If $x\in (G/P)^T$, then the closure of 
the $B$-orbit $Bx$ is called the \SV\ in $G/P$ associated to $x$. 
This \SV\ will be denoted throughout by $X(x)$. We will use the well known fact that
the \tfp s in $G/B$ are in one to one correspondence with the 
elements of the Weyl group $W=N_G(T)/T$, so we don't distinguish
between elements of $W$ and fixed points in $G/B$.

\SVS\  are in general singular, 
and it's an old problem, inspired by a classical paper \cite{chev} of Chevalley, 
to describe their singular loci (or, equivalently, their
smooth points). A related problem, with interesting consequences in  
representation theory, is to determine the locus of rationally smooth
points of a \SV\ (cf. \cite{kl}). In fact, if $G$ is defined over $\C$ and simply laced
(i.e. every simple factor is of type $A$, $D$ or $E$), then all
rationally smooth points of any \SV\ in $G/P$ are in fact smooth (see \cite{ck}).

In this paper, we  consider the 
singular locus of a \SV\ in an arbitrary 
$G/P$, where $G$ does not contain any $G_2$-factors.
Our results are an outgrowth of \cite{ck}, where we used \pt s (defined below) to
characterize the \tfp s in the singular locus of an irreducible $T$-stable  
subvariety $X\subset G/P$ (a $T$-variety for short). 

Since every $B$-orbit meets $(G/B)^T$,  the singular locus of a \SV \ $X$ 
consists of the $B$-orbits of its singular  \tfp s.
Here we can make use of the well known
natural ordering order on $W=(G/B)^T$: if  $x,y\in W$,
then $x<y $ \iff $X(x)\subset X(y)$ but $x\ne y$. This ordering
can also be described as the Bruhat-Chevalley order on $W$ (cf. \cite{chev}).
It has two nice properties. First, if $X=X(w)$, then
$X^T$ is the interval $[e,w]=\{x\in W \mid x\le w\}$. Secondly, 
if $X$ is smooth (resp. singular) at $y\in X^T$, 
then it is smooth (resp. singular) at all
$z\in X^T$ with $z>y$ (resp. $z\le y$). Hence the problem
of determining the singular locus of $X$ boils down 
to identifying the maximal singular elements $y$ of $[e,w]$. 
Such a $y$ is called a {\em maximal singularity} of $X$.


Let us now recall the notion of a \pt\ from \cite{ck}. 
Suppose for the moment that $X$ is a $T$-variety in $G/P$,
and let $E(X,x)$ denote the set of
$T$-stable curves in $X$ which contain the point $x\in X^T$. 
An element $C$ of $E(X,x)$ is called {\em good} if $C$
contains a smooth point of $X$. The {\em \pt} of
$X$ along a good $C\in E(X,x)$ can be defined as
$$\tau_C(X,x)=\lim _{z\to x}T_z(X)\quad(z\in C \backslash C^T),$$
where $T_z(X)$ is the Zariski tangent space of $X$ at $z$.
Each \pt\ $\tau_C(X,x)$ is a $T$-stable subspace of $T_x(X)$
\st $\dim \tau_C(X,x) =\dim X$. 
Clearly, $X$ is smooth at $x$ \iff $\tau_C(X,x)=T_x(X)$. 
Each $T_x(C)$ is a $T$-stable line in $T_x(X)$, and the $T$-weights of these lines are 
certain elements of the root system $\Phi$ of the pair $(G,T)$. 
Furthermore, if $T_x(C)$ has weight $\A\in \Phi$, then $C^T=
\{x,r_\A x\}$, where $r_\A \in W$ is the reflection corresponding to $\A$.
A $T$-curve $C\in E(X,x)$ is called 
{\em short} or {\em long} according to whether the $T$-weight of 
its tangent line $T_x(X)$ is short or long in $\Phi$.

Now put $TE(X,x)=\sum_{C\in E(X,x)} T_x(C)$, and  let
$\tx$ be the linear span of the reduced tangent cone of $X$ at $x$. 
Clearly, $TE(X,x)\subset \tx$. In fact, $TE(X,x)= \tx$ if $G$
is simply laced and $X$ is a \SV\ \cite{cp}. 
(We will extend this to arbitrary $T$-varieties in $G/P$
in \S 3.)  A key question is whether $\tau_C(X,x)\subset TE(X,x)$.
The following result combining Theorems 1.4 and 1.6 of \cite{ck}
gives a partial answer in this setting. 

\begin{theorem}\label{TH1}  Assume $G$ has no $G_2$-factors, and let $X$ be a $T$-variety in $G/P$ 
which is \cm \ at the point $x\in X^T$.
Then $X$ is smooth at $x$ \iff $\tau_C(X,x)=TE(X,x)$
for at least  one good $T$-curve $C$. Moreover, 
\begin{equation}\label{PTINTAU} \tau_C(X,x)\subset \tx
\end{equation}
for every good $C\in E(X,x)$. In fact, if $C$ is good and also short, 
then $\tau_C(X,x)\subset TE(X,x)$.
\end{theorem}

\begin{remark} In fact, the first assertion is true even if $G$ has 
$G_2$ factors. However, as we show in Example \ref{G2eg},
the identity (\ref{PTINTAU}) fails in $G_2/B$, as does the assertion
about short curves. 
\end{remark}

Returning to the case of a \SV\ $X$, this result reduces the problem problem of classifying the 
singular locus of $X$ to 
describing  $\tx$ when $E(X,x)$ contains a good  curve
and $G$ doesn't have $G_2$ factors. Our main goal is to solve this
for  a maximal singularity of $X$. 


The first of our two characterizations of $\tx$
at a maximal singularity goes as follows.

\begin{theorem} \label{SUM_TAU} Assume $G$ has no $G_2$-factors, and suppose 
$X$ is a \SV \ in $G/P$.
Then, for any  $x\in X^T$ which is either smooth in $X$ or a maximal 
singularity of $X$, 
\begin{equation}\label{THETA=SUM_TAU} \Theta_x(X)=\sum_{C}\tau_C(X,x),
\end{equation}
where the sum is over all $C\in E(X,x)$ with $C^T =\{x,y\}$, where $y>x$. 
\end{theorem}

The proof is given in \S \ref{STC}. We also include there a counter-example 
(Example \ref{G2eg}) to 
Theorem \ref{SUM_TAU} in the $G_2/B$ case. 
Note that if $x$ is maximal (or smooth), all $T$-curves $C$ \st 
$C^T =\{x,y\}$ and $y>x$ are good.
An algorithm for computing $\tau_C(X,x)$ was given in \cite[\S 8]{ck}, so  
(\ref{THETA=SUM_TAU}) gives an explicit method of computing $\tx$.

Our second characterization uses Theorem \ref{SUM_TAU} to
get an expression for $\tx$ in which  \pt s are out of the
picture.  For this, we need to bring in some more concepts.
For any $x\in X^T$, let $B_x \subset B$ be the isotropy subgroup of $x$. 
That is, $B_x$ is the subgroup of $B$ namely the subgroup
generated by $T$ and all root subgroups $U_\A$ of $B$ which fix $x$.
As usual, a root $\A$ \st $U_\A \subset B$ is taken 
to be positive, and we write $\A >0$. The condition that $U_\A x=x$ is equivalent to  
$x^{-1}(\A)>0$. Thus, for any \SV\ $X$ in $G/P$, 
$\tx$ is a $B_x$-submodule of $T_x(G/P)$. 
The {\em isotropy submodule} of $\tx$ is the smallest $B_x$-submodule 
$\TX$ of $T_x(X)$ which contains $TE(X,x)$. 

We will show  that if $C\in E(X,x)$ is good, then the roots which correspond
to $T$-lines in the $T$-module $\tau_C(X,x)/(\TX\cap \tau_C(X,x))$ can be explicitly
described in terms of the notion of an orthogonal
$B_2$-pair, which is now defined. 
For each $\g \in \Phi$, let $\fg_\g$ 
denote the $T$-stable line in the Lie algebra $\fg=\mathrm{Lie}(G)$ 
of weight $\g.$ In other words, $\fg_\g$ is the root line of weight $\g$.

\begin{definition}\label{B2ROOT} Let $X=X(w)$ be a \SV \ in $G/B$, and assume $x< w$.
Suppose $\mu$ and $\phi$ are long, positive orthogonal roots \st
the following three conditions hold:
\begin{itemize}

\item[$(i)$]  $\fg_{-\mu} \oplus \fg_{-\phi} \subset  TE(X(w),x)$ (hence 
$x<r_\mu x, ~r_\phi x \le w$),

\medskip
\item[$(ii)$]  there exists a subroot
system $\Phi'$ of $\Phi$ of type $B_2$  containing $\mu$ and $\phi$, and

\medskip
\item[$(iii)$]
if $\A$ and $\B$ form the unique basis of $\Phi'$ contained in
$\Phi^+ \cap \Phi'$
with $\A$ short and $\B$ long, then
\begin{equation*}
r_\A x<x, \quad \mathrm{and} \quad r_\A r_\B x \le w.
\end{equation*}
\end{itemize}
Then we say that $\{\mu,\phi\}$ form an {\em orthogonal $B_2$-pair for 
$X$ at} $x.$
\end{definition}

Our second characterization of $\tx$ at a maximal singularity goes as follows.
\begin{theorem}\label{THB2PAIR} Assume $G$ has no $G_2$-factors,
and suppose $x$ is a maximal singularity of a \SV\ $X$ in $G/B.$ 
Then for each $T$-weight $\g$ of the quotient $\tx/\TX$, there exists 
an orthogonal $B_2$-pair $\{\mu,\phi\}$ for $X$ at $x$ \st
\begin{equation}\label{AVEWEIGHT} \g=-1/2(\mu +\phi).
\end{equation}
In other words, at a maximal singularity of $X$, every $T$-weight of $\tx$
not in $\TX$ is a weight arising from a $B_2$-pair at $x$ as in $(\ref{AVEWEIGHT})$.
\end{theorem}

This is proved in \S \ref{PTSV}.
We also obtain a necessary and sufficient condition for 
a \tfp\ $x$ of a \SV\ to be a smooth point, which is also proved in \S \ref{PTSV}.

\begin{theorem}\label{BGALG} Assume $G$ has no $G_2$-factors, let $X$ be a \SV \ in $G/B$, 
and let $x\in X^T$. Then $X$ is smooth at $x$ \iff the following
three conditions hold.

\ms
$(1)$ $|E(X,x)|=\dim X$, and some $C\in E(X,x)$ is good.

\ms
$(2)$ We have $\TX =TE(X,x)$, and

\ms
$(3)$ If  $\{\mu,\phi\}$ is an orthogonal $B_2$-pair  for $X$ at $x$,
and  $\g=-1/2(\mu +\phi)$, then $\fg_\g \subset TE(X,x)$.
Consequently, $r_\g x \le w$.
\end{theorem}

\begin{corollary} \label{ALGCOR} There exists a 
non-recursive algorithm involving only the  \bg\  and the root system $\Phi$
which classifies the smooth \tfp s of a \SV\ in $G/B$. 
\end{corollary}

The notion of an orthogonal $B_2$-pair arises from 
the \SV \ $X=X(r_\A r_\B r_\A)$ in $B_2 /B$,
where $\A$ and $\B$ are the short and long simple roots in $B_2$.
The \tfp\ $x=r_\A$ is the unique maximal singularity of $X$.
Now the weights of $TE(X,x)$ are $\A, -\B$ and $-(\B +2\A)$.
Furthermore, $B_x$ is generated by $T$, $U_\B$, $U_{\A+\B}$ and $U_{2\A+\B}$,
so it is easy to see that $TE(X,x)$ is already a $B_x$-submodule of $T_x(X)$.
But $\{\B, \B+2\A\}$ give an orthogonal $B_2$-pair at $x$
\st $\fg_\g\subset \tx/TE(X,x)$, where $\g=-1/2(\mu +\phi)=-(\A +\B)$.
(See Example \ref{B2} and \cite{ck} for more details.)

\bigskip

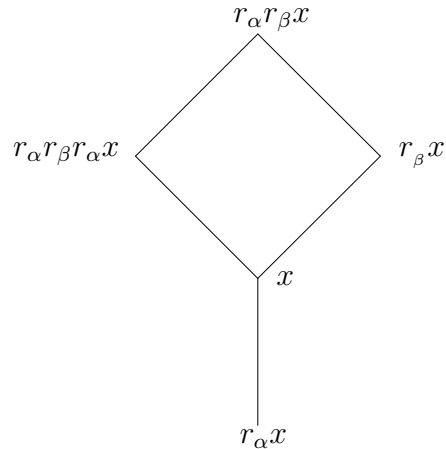
\begin{center}
\setlength{\unitlength}{0.00043333in}
\begingroup\makeatletter\ifx\SetFigFont\undefined%
\gdef\SetFigFont#1#2#3#4#5{%
  \reset@font\fontsize{#1}{#2pt}%
  \fontfamily{#3}\fontseries{#4}\fontshape{#5}%
  \selectfont}%
\fi\endgroup%
{\renewcommand{\dashlinestretch}{30}
\begin{picture}(5186,5370)(0,-10)
\path(3000,2025)(4500,3525)
\path(4500,3525)(3000,5025)
\path(3000,2025)(1500,3525)
\path(1500,3525)(3000,5025)
\path(3000,2025)(3000,225)
\put(3225,1950){\makebox(0,0)[lb]{\smash{{{\SetFigFont{12}{14.4}{\rmdefault}{\mddefault}{\updefault}$x$}}}}}
\put(4725,3525){\makebox(0,0)[lb]{\smash{{{\SetFigFont{12}{14.4}{\rmdefault}{\mddefault}{\updefault}$r_{_\B} x$}}}}}
\put(2700,5175){\makebox(0,0)[lb]{\smash{{{\SetFigFont{12}{14.4}{\rmdefault}{\mddefault}{\updefault}$r_\A r_\B x$}}}}}
\put(0,3525){\makebox(0,0)[lb]{\smash{{{\SetFigFont{12}{14.4}{\rmdefault}{\mddefault}{\updefault}$r_{\A} r_{\B} r_{\A} x$}}}}}
\put(2775,0){\makebox(0,0)[lb]{\smash{{{\SetFigFont{12}{14.4}{\rmdefault}{\mddefault}{\updefault}$r_{\A} x$}}}}}
\end{picture}
}

\end{center}

\medskip
\begin{center} {\bf Figure 1:} $\A$ and $\B$ are the short and long\\  simple 
roots in a $\Phi^+(B_2)$ containing $\{\mu,\phi\}$.
\end{center}

Figure 1 illustrates the portion of \bg\ of a \SV\
$X$ arising from an orthogonal $B_2$-pair $\{\mu,\phi\}$ at $x.$ 
If $x$ is on a good $T$-curve and there is no edge in $\G(X)$ 
at $x$ corresponding to a $T$-curve $C$ with 
$x<r_\g x \le w$, where $\g=-1/2(\mu +\phi)$, then $X$ is singular at $x$.

Let us describe the algorithm of Corollary \ref{ALGCOR}.
Suppose we want to determine whether a \SV\ $X=X(w)$ is smooth at $x\in X^T.$ 
Consider any  descending path
$$w > x_1 > x_2 > \cdots > x_m >x$$
in $\G(X).$ If $X$ is singular at any $x_i$, then it is singular at $x$.
Thus, suppose  $X$ is smooth at $x_{m}$. Then the edge
$x_{m} x$  is a good $T$-curve in $X$,  
so it suffices to check the conditions of  Theorem \ref{BGALG} for this $T$-curve. 
Checking that the number of $T$-curves is $\dim X$ amounts  to showing $|\{\g >0 \mid r_\g x\le w\}|=\ell (w),$
where $\ell(w)$ is the length of $w$, since $\ell(w)=\dim X(w)$.
Verifying the second condition amounts to
showing that $TE(X,x)$ is $B_x$-stable.
This requires verifying that if $\fg_\g \subset TE(X,x)$, then $\fg_{\g 
+\A} \subset TE(X,x)$
for all $\A >0$ \st $x^{-1}(\A)>0$, $\g +\A \in \Phi$ and $x^{-1}(\A +\g)<0$. 
The third condition is verifiable from the \bg\ at $x$,
so the algorithm involves only
$\Phi$ and $\G(X)$. The non-recursivity is due to the fact that
we only need to consider a single path in $\G(X)$ from $w$ to $x$.


It might also be useful to remark that unlike checking whether a \SV\ $X$
is smooth at a fixed point $x$ , checking for rational smoothness at $x$ via 
the \bg\ requires that one count the number of edges
in $\G(X)$ at all vertices $y\ge x$  \cite{cp}.
Therefore it appears to be easier to use the \bg\ to identify the smooth points
than the rationally smooth points. 
B. Boe and W. Graham
have formulated the following lookup conjecture: 
a \SV\ $X$ in $G/P$  is rationally smooth
at $x$ \iff $|E(X,y)|=\dim X$ for all $y$ on an edge of $\G(X)$
containing $x$.
Some special cases of the lookup conjecture
are verified in \cite{bg}, but the general conjecture
is open. Theorem \ref{BGALG}  says that as far as smoothness
is concerned, one has to examine $\G(X)$ two steps above and one step 
below $x$. This might be considered somewhat unexpected.

Finally, let us mention that this paper has connections with the work
of S. Billey and A. Postnikov \cite{bp} and very likely also 
S. Billey and T. Braden \cite{bb}. However, unlike the situation in 
\cite{bp}, our results do not say anything in the 
$G_2$ case, as noted in Remark \ref{G2eg}.

\section{Preliminaries}\label{prelim}
We will throughout use the terminolgy and notation of 
\cite{ck}, some of which was already introduced in \S \ref{intro}. 
In particular, the $G_2$-hypothesis is always in effect. 

Let us first recall  some of the standard facts and notations 
concerning roots, weights, $T$-curves and so on. 
The $T$-fixed point set of a $T$-variety  $X\subset G/P$  
is denoted by $X^T$. It's well known that the 
mapping $w\to n_wB$ is a bijection the Weyl group $W=N_G(T)/T$ 
of $(G,T)$ with $(G/B)^T$, so we assume $W=(G/B)^T$.
The projection $\pi:G/B \to G/P$ is an equivariant closed morphism,
so $(G/P)^T$ may be identified with $W/W_P$,  
$W_P$ being the parabolic subgroup of $W$ associated to $P$.
The elements of $W/W_P$ thus parameterize the \SVS\ in $G/P$. 

Every $T$-curve in $E(X,x)$ has the form $C=\overline{U_\A x}$ 
for a unique root $\A \in \Phi$. If $P=B$, then $C^T=\{x, r_\A x\}$.
If $X$ is a \SV \ in $G/B$, say $X=X(w)$,
then $C=\overline{U_\A x} \subset X$ \iff both $x,r_\A x \le w$. By 
\cite[LEMMA A]{cp}, $|E(X,x)|\ge \dim X$ for every $T$-variety $X$.
Furthermore, every $T$-curve in $G/P$ is the image of a $T$-curve in $G/B$
under the closed morphism $\pi:G/B \to G/P$. Also,
recall that as $T$-modules,
$$T_x(G/B) =\bigoplus_{x^{-1}(\g)<0} \fg_\g.$$

A property of $T$-varieties in $G/P$, used freely
throughout the paper is the following: each $T$-fixed point $x \in
G/P$ is \emph{attractive}; that is, all the weights of the tangent
space $T_x(G/P)$ lie on one side of a hyperplane in $X(T)$, and in
addition, each fixed point $x$ has a $T$-stable open affine
neighborhood. Since $X$ is irreducible
and any $x \in X^T$ is attractive, the affine open $T$-stable
neighborhood of $x$ is unique. It will be denoted by $X_x$. It is
well known, and not hard to see, that there is a closed
$T$-equivariant embedding of $X_x$ into the tangent space $T_x(X)$
of $X$ at $x$, thanks again to the fact that $x$ is attractive.

Assuming $X_x \subset T_x(X)$, it follows that, for any
$T$-stable line $L \subset T_x(X)$, we may choose a linear
equivariant projection $T_x(X) \RA L$ and restrict it to $X_x$.
Identifying $L$ with $\mathbb A^1_k$ we thus obtain a regular
function $f \in k[X_x]$, which is a $T$-eigenvector of weight $-\alpha$ if
$L$ has weight $\alpha$. We say $f$ \emph{corresponds to $L$}, if
it is obtained in the described way.

\section{Some General Results on $\tx$}\label{SL}

In this section,  we will establish some general properties of an
arbitrary $T$-variety $X$ in $G/P$. For Schubert varieties these
properties are well known (see \cite{car}).  Let $\T _x(X)$ be the
reduced tangent cone to $X$ at any $x\in X^T$, so
$\tx=\mathrm{span}_k (\T_x(X))$.

\begin{theorem} \label{TC=TL} Suppose $G$ has no $G_2$-factors. Let
  $L=\fg_\omega
\subset \Theta_x(X)$ be a $T$-stable line with weight $\omega$.
Then the following hold.

\ms

\begin{itemize}

\ms

\item[$(i)$] If $\omega$ is long, then $L\subset TE(X,x)$. Otherwise,
there exist  roots $\A,\B$ \st
$\fg_\A,\fg_\B \subset TE(X,x)$ and
\begin{equation}\label{AVE}
\omega = \frac{1}{2}(\A + \B).
\end{equation}



\item[$(ii)$] In particular, if $G$ is simply laced, then $\tx =TE(X,x)$.

\item[$(iii)$] If $X$ is a \SV \ and $L$ does not correspond
to a $T$-curve, then $\A$ and $\B$ are long negative orthogonal roots
in a copy of $B_2 \subset \Phi$.

\end{itemize}

\end{theorem}

\begin{proof}
Let $z \in k[X_x]$ be a $T$-eigenfunction corresponding to $L$ and
let $x_1,x_2,\dots,x_n \in k[X_x]$ be $T$-eigenfunctions which
correspond to the $T$-curves $C_1,C_2, \dots, C_n$ through $x$.
Notice that since $X_x \subset T_x(X)$ each $T$-curve $C \in
E(X_x,x)$ is in fact a coordinate line in $T_x(X)$. This follows
from the fact that all $T$-curves are smooth and no two
$T$-weights of $T_x(X)$ are proportional. Let $\tilde x_i$,
resp. $\tilde z$ denote linear projections $T_x(X) \RA T_x(C_i)$, resp.
$T_x(X) \RA L$, which restrict to $x_i, z \in k[X_x]$.

Since the (restriction of the) projection $X_x \RA \bigoplus_C
T_x(C) = TE(X,x)$ has a finite fibre over $0$, $k[X_x]$ is a
finite $k[x_1,x_2, \dots , x_n]$-module by the graded version of
Nakayama's Lemma. In particular $z \in k[X_x]$ is integral over
$k[x_1,\dots,x_n]$, so we obtain a relation
\begin{equation} \label{INTEQ}
z^ N = p_{N-1} z^{N-1} + p_{N-2}z^{N-2} + \dots + p_1 z + p_0,
\end{equation}
where $N$ is a suitable integer and $p_i\in k[x_1,\dots,x_n]$.
Without loss of generality we may assume that every summand on the
right hand side is a $T$-eigenvector with weight $N \omega$. Let
$P_i \in k[\tilde x_1, \dots, \tilde x_n]$ be a polynomial
restricting to $p_i$, having the same weight $(N-i)\omega$ as
$p_i$. Then every monomial $m$ of $P_i$ has this weight too. If
for all $i$ every such monomial $m$ has degree $\deg m > N - i$,
then $p_i z^{N- i}$ is an element of $M^{N + 1}$, where $M$ is the
maximal ideal of $x$ in $k[X_x]$. This means that $\tilde z$
vanishes on the tangent cone of $X_x$, so $L \not \subset
\Theta_x(X)$, which is a contradiction.
Thus, there is an $i$ and a monomial $m$ of $P_i$, such that
$\deg m \leq d = N - i$. Let $m = c
\tilde x_1^{d_1} \tilde x_2 ^{d_2}
\dots \tilde x_n^{d_n}$, with integers $d_j$ and a nonzero $c \in
k$. So $\sum_j d_j \leq d$. Let $\A_j$ be the weight of $\tilde
x_j$. Then we have
\begin{equation*}
  d \omega = \sum d_j \A_j
\end{equation*}
After choosing a new index, if necessary, we may assume that $d_j
\not = 0$ for all $j$. Let $(~,~)$ be a Killing form on
$X(T)\otimes {\mathbb R}$ which induces the length function on
$\Phi$. We have to consider two cases. First suppose that $\omega$
is a long root, with length say $l$. Then $(\A_j, \omega) \leq
l^2$ with equality if and only if $\A_j = \omega$. Thus, $d l^2 =
\sum d_j (\A_j,\omega) \leq d \max_{j} (\A_j, \omega) \leq d l^2$
and so there is a $j$ with $\A_j = \omega$ and we are done, since
this implies $\tilde z = \tilde x_j$. Hence, $L=C_j$.

Now suppose $\omega$ is short, with its length also denoted $l$.
In this case $(\A_j, \omega) \leq l^2$. Since $dl^2 =
d(\omega,\omega) = \sum_j d_j (\A_j, \omega)$ and since $\sum d_j
\leq d$, it follows that all $\A_j$ satisfy $(\A_j, \omega) =
l^2$. If there is a $j$ such that $\A_j = \omega$, then, as above,
we are done. Otherwise for each $j$, $\A_j$ is long, and $\A_j$
and $\omega$ are contained in a copy $B(j) \subset \Phi$ of $B_2$.
There is a long root $\B_j \in B(j)$ with $\A_j + \B_j = 2\omega$.
We have to show that there are $j_0$ and $j_1$ so that $\B_{j_0} =
\A_{j_1}$. Fix $j_0 = 1$ and let $\A = \A_1$, $\B = \B_1$. Then
$(\A,\B) = 0$.
  This gives us the result: $d l^2 = d(\omega, \B) = 0 + \sum_{j>1}
(\A_j, \B)$. Now if all $(\A_j,\B)$ are less or equal $l^2$, this
last equation cannot hold, since $\sum_{j>1}d_j < M$. We conclude
that there is a $j_1$ so that $(\A_{j_1}, \B) = 2l^2$ (the squared
long root length), hence $\A_{j_1} = \B$, and we are through with
$(i)$.

The proof of $(ii)$ is obvious. For $(iii)$, let $S$ be the slice
(cf. \cite[]{ck}) to $X(w)$ at $x$. Then, locally, $X=S\times Bx$,
where the weights of $TE(S,x)$ consist of the roots $\A<0$ \st
$x<r_\A x \le w$. Since $L\not \subset TE(X,x)$, the only
possibility is that $L \subset \Theta_x(S)$ because $Bx$ is smooth
(and so $TE(Bx,x) = \Theta_x(Bx)$) and $\Theta_x(X) = \Theta_x(S)
\oplus \Theta_x(Bx)$. No we may apply part $(i)$ to $S$.
\end{proof}

The following generalizes a well known property of Schubert varieties to
  arbitrary $T$-varieties.



\begin{corollary}\label{TC=TLCOR}
Suppose  $L$ is a $T$-invariant line $\T_x(X)$. Then $L\subset TE(X,x)$.
\end{corollary}

\begin{proof} We have already shown that in equation (\ref{INTEQ}), some
$P_i$ contains a monomial of degree at most $d=i$.
Taking homogeneous parts of degree $N$ in  (\ref{INTEQ}), we therefore get
a homogeneous polynomial
$$f=\tilde{z}^N-\sum P_j \tilde{z}^{N_j}.
$$
vanishing on $\T_x(X)$. Hence $f(L)=0$. But as  $\tilde{z}(L)\ne
0$, this implies some $P_j(L)\ne 0$ as well, which means that
$\tilde z$ occurs in a monomial of $P_j$, hence $L\subset TE(X,x)$
by the construction of the $P_j$.
\end{proof}


An interesting consequence of Corollary \ref{TC=TLCOR}
is that the linear spans of
the tangent cones of two $T$-varieties behave nicely under
intersections.
\begin{corollary} \label{INT} Suppose the $G$ is simply laced and that
$X$ and $Y$ are $T$-varietes in $G/P$. Suppose also that $x\in
(X\cap Y)^T $. Then
$$\Theta_x(X\cap Y)=\Theta_x(X)\cap \Theta_x(Y).$$
Consequently, if both
$X$ and $Y$ are nonsingular at $x$, then
$X\cap Y$ is nonsingular at $x$ \iff $|E(X\cap Y,x)|=\dim (X\cap Y)$.
\end{corollary}

\begin{proof} The first claim is clear since $E(X,x)\cap E(Y,x)=E(X\cap 
Y, x)$.
For the second, note that if $X$ and $Y$ are nonsingular at $x$, then
\begin{eqnarray*}
T_x(X)\cap T_x(Y)&=&\Theta_x(X)\cap \Theta_x(Y)\\
                   &=&\Theta_x(X\cap Y)\\
                   &\subset & T_x(X\cap Y)\\
                   &\subset & T_x(X)\cap T_x(Y)
\end{eqnarray*}
Hence $\dim T_x(X\cap Y) =|E(X\cap Y)|,$ and the result follows.
\end{proof}

For example, it follows that in the simply laced setting, the
intersection of a Schubert variety $X(w)$ and a dual Schubert
variety $Y(v)=\overline{B^-y}$ is nonsingular at any $x\in [v,w]$
as long as $X(w)$ and $Y(v)$  are each nonsingular at $x$.

\section{$\tx$
at a Maximal Singularity} \label{STC}

The aim of this section is to prove Theorem \ref{SUM_TAU}.
In fact,  we will derive it as a consequence
of a general result about the connection between
$\tau_C(X,x)$ and $\Theta_x(X)$ for an arbitrary $T$-variety
in $G/P$ assuming  $x$ is at worst
an isolated singularity.

\begin{theorem} \label{THETA_SUB_TAU}
Suppose $X \subset G/P$ is a $T$-variety, where $G$ has no
$G_2$-factors. Then for each $x \in X^T$, we have
$$\Theta_x(X) \subset \tau(X,x):=\sum _{C\in E(X,x)} \tau_C(X,x).$$
In particular, if $x$ is either smooth in $X$ or  an isolated singularity,
then
$$\Theta_x(X)=\sum _{C\in E(X,x)} \tau_C(X,x).$$
\end{theorem}

Before proving Theorem \ref{THETA_SUB_TAU}, we will give
the proof of Theorem \ref{SUM_TAU}.
\begin{proof}[Proof of Theorem \ref{SUM_TAU}] The result is
obvious if $x$ is smooth, so assume it is a maximal singularity.
Then there exists a slice representation
$X_x=S\times Bx$, where $S$ has an isolated singularity at $x$
and  $E(S,x)$ consists of the $T$-curves in $X$ containing a smooth point
of $X_x$.
To get the result, we apply Theorem \ref{THETA_SUB_TAU}
to $S$ and use the fact that $\tx =\Theta_x(S)\oplus \Theta_x(Bx)$.
Indeed,
$$\Theta_x(S)\oplus \Theta_x(Bx)=\sum_{C\in E(S,x)} \tau_C(S,x) \oplus
TE(Bx,x),$$
so it suffices to show that $TE(Bx,x)\subset  \tau_C(X,x)$ for any ${C\in
E(S,x)}$
since clearly  $\tau_C(S,x)\subset \tau_C(X,x)$.
Let $\fg_\g \subset TE(Bx,x)$. Then there is a curve $D \subset
Bx$ with $\fg_\g = T_x(D)$. In fact, $D = U_\g x$. Thus, the
smooth $T$-stable surface $\Sigma = C \times D \subset X_x =
S\times Bx$, and Proposition 3.4 of \cite{ck} implies $\fg_\g
\subset \tau_C(\Sigma,x) \subset \tau_C(X,x)$.



$\Sigma=\overline{U_\g C }$ is a $T$-surface

hard to see




\end{proof}

The proof of Theorem  \ref{THETA_SUB_TAU} will use
several lemmas. To begin with,
let $R$ be a Noetherian graded commutative ring, with irrelevant
ideal $I = \bigoplus_{d > 0}  R_d$. Then $\bigcap_{l> 0} I^l = 0$.
Thus, for each $r \in R \setminus \{0\}$ there is an $l > 0$ such
that $r \in I^l \setminus I^{l+1}$. We set $\IN(r) = r + I^{l+1}
\in I^l/I^{l+1} \subset \gr R = \gr_I R$, and $\IN(0) = 0 \in \gr
R$. Recall that for $r,s \in R$, $\IN(r)\IN(s) = \IN(rs)$ or
$\IN(r)\IN(s) = 0$. We say $r \in R$ \emph{vanishes on the tangent
cone} if $\IN(r)$ does, i.e. if $\IN(r)$ is nilpotent. In the case
that $R$ is the coordinate ring of an affine variety $Z$ with
regular $\mathbb G_m$-action such that $I$ corresponds to a
maximal ideal and hence to an attractive $\mathbb G_m$-fixed point
$z$, then $\IN(r)$ induces indeed a function on the reduced
tangent cone of $Z$ at $z$, and $r$ vanishes on the tangent cone
if and only if this function does. In what follows we will
consider closed and $T$-stable subvarieties of $T_x(X)$. We
therefore choose a one parameter subgroup $\lambda$ of $T$, such
that $\lim_{t\RA 0}\lambda(t)v = 0$ for all $v \in T_x(X)$. Then
the $\mathbb G_m$-action by $\lambda^{-1}$ induces a (positive)
grading of $k[T_x(X)]$ which carries over to any $T$-stable closed
subvariety (note that the grading induced by $\lambda$ would be
negative).

For convenience we extend the definition of $\Theta_x(Z)$ also to
reducible varieties. Also notice that $\Theta_x(Z)$ may be
canonically identified with $T_0(\mathfrak T_z(Z)) \subset
T_z(Z)$. To set up an induction on the dimension of $X$, we need
the following

\begin{lemma} \label{THETA_SUM}Let $Z \subset T_x(X)$ be a closed
$T$-stable
subvariety with $Z = Z_1 \cup Z_2 \cup \dots \cup Z_d$ the decomposition
into irreducible comoponents.
  Then $$\Theta_0(Z) = \Theta_0(Z_1) + \Theta_0(Z_2) + \dots
  + \Theta_0(Z_d).$$
\end{lemma}

\begin{proof}
Since every component $Z_i$ of $Z$ is $T$-stable it has to contain $0$.
Therefore the proof is a simple consequence of the following well known
fact:
if a variety $Y = A \cup B$ is the union of two closed subvarieties then
for
every point $x$ in the intersection $A \cap B$ we have $\T_{x}(Y) = \T_x
(A)
\cup \T_x (B)$.
\end{proof}

Let $Z \subset T_x(X)$ be an irreducible $T$-stable subvariety,
and let $L \subset \Theta_0(Z)$ be a $T$-stable line with weight
$\omega$, say. Moreover, suppose $\omega$ is short with respect to
a Killing form $(~,~)$ on $X(T)$. Denote by $z \in k[Z]$ the
restriction of a linear $T$-equivariant projection $T_x(X)
\rightarrow L \cong \mathbb A^1_k$. We fix $z$ for the moment.

\begin{lemma}With the preceding notation, let $f \in k[Z]$ correspond to
another $T$-equivariant linear projection onto some
line $L^\prime \subset T_x(X)$. Then $z$ vanishes on the tangent
cone of $\mathcal V(f)$ if and only if $\IN(z)^l = \IN(h)\IN(f)$
for some positive integer $l$ and a suitable $T$-eigenvector $h
\in k[Z]$.
\end{lemma}

\begin{proof}
The if is clear, so suppose $z$ vanishes on the tangent cone of
$\mathcal V(f)$. By definition this means that there is an integer
$l$ and there are elements $g_1,g_2,\dots, g_r \in I(\mathcal
V(f))$, the ideal of $\mathcal V(f)$, such that $\IN(z)^l =
a_1\IN(g_1) + a_2 \IN(g_2) +\dots + a_r\IN(g_r)$ for suitable $a_i
\in \gr k[Z]$. Since $\IN(z)$ is homogeneous and since
$\IN(I(\mathcal V(f))$ is an homogeneous ideal, we may assume that
all of the $a_i$ are homogeneous as well. Moreover the $a_i$ and
$g_i$ may be chosen to be $T$-eigenvectors. Omitting any indices
$i$ for which $a_i\IN(g_i) = 0$ we may lift the $a_i$
equivariantly to $\bar a_i \in k[Z]$, such that $\IN(\bar a_i) =
a_i$. Then we have $0 \neq \IN(\bar a_i)\IN(g_i) = \IN(a_ig_i)$.
Leaving out degrees different from $l$ we may assume that $\sum
\IN(\bar a_i)\IN(g_i) = \IN(\sum \bar a_i g_i)$. Now $\sum \bar
a_i g_i$ is a $T$-eigenvector $g$ contained in the ideal of
$\mathcal V(f)$. A suitable $n$th power of $g$ is contained in
$fk[Z]$. $\IN(z)$ is not nilpotent, and due to $\IN(z)^l = \IN(g)$
  also $\IN(g)$ is not nilpotent, therefore $\IN(g)^n = \IN(g^n)$.
Replacing
$l$ by $nl$ we may
  assume that $\IN(z)^l = \IN(g)$ for a $g \in fk[Z]$. In other
  words $\IN(z)^l = \IN(hf)$ for a suitable $T$-eigenvector $h \in
  k[Z]$.

  It remains to show that $\IN(hf) = \IN(h)\IN(f)$ which is
  equivalent to $\IN(h)\IN(f) \neq 0$. So suppose that
  $\IN(h)\IN(f) = 0$. This means that $h \not \in M^{l-1}$ with $M$ the
maximal ideal of $0$.
  Otherwise $\IN(h)\IN(f)$ would equal $\IN(hf)$ by definition. We
  conclude that $h \in M^n$ for some $n < l-1$, implying that there is a
  a homogeneous polynomial $P$ in some linear $T$-homogeneous coordinates
$x_1,x_2,\dots,x_m$ of
  $T_x(X)$ of the same $T$-weight as $h$, and of degree $n$, such
  that, restricted to $Z$, $h = P$ modulo $M^{n+1}$. By the definition of
$f$
we may even assume that $x_1$ restricted to $Z$ is $f$. Replacing $P$ by
any
  monomial of $P$ and letting $d_i$ be the degree of $x_i$ in $P$,
  we see that $l\omega = \alpha_1+ \sum d_i \alpha_i$ with $\alpha_i$ the
  weight of $x_i$. Applying $(\cdot,\omega)$ on both sides this
  gives $l(\omega,\omega) = (\alpha_1,\omega)+ \sum d_i
  (\alpha_i,\omega)$. Since $(\alpha_1,\omega) \leq
  (\omega,\omega)$ for all $i$ this is impossible since $n = \sum
  d_i < l -1$. Hence the claim.
\end{proof}

As an easy consequence we get

\begin{lemma}
If $Z$ and $z$ are as above, and $f$ corresponds to the projection
to any other $T$-stable line of $T_x(X)$ with a short weight, then
$z$ does not vanish on the tangent cone of $\mathcal V(f)$.
\end{lemma}

\begin{proof}
By the last Lemma we know, that if $z$ vanishes on the tangent
cone of $\mathcal V(f)$, there is a $T$-eigenvector $h \in k[Z]$
such that $\IN(z) = \IN(h)^l\IN(f)$. Choosing a monomial as in the
proof of the previous Lemma, we get a relation $l\omega = \alpha_1
+ \sum d_i \alpha_i$ with $\sum d_i = l-1$. But $(\alpha_1,\omega)
< (\omega,\omega)$, because $\alpha_1$ is short, and
$(\alpha_i,\omega) \leq (\omega,\omega)$ for all $i$, so no such
relation exists.
\end{proof}

For reasons which will become clear in the proof of the Theorem,
we now restrict our attention to varieties $Z$ in $T_x(X)$ such
that $T_0(Z)$ contains exactly one $T$-stable line with a short
weight.

\begin{lemma}
If $L \subset \Theta_0(Z)$ is the only line in $T_0(Z)$, and if $C
\in E(Z,0)$ is any $T$-curve, then $L \subset T_p(Z)$ for all $p
\in C^o = C \setminus \{0\}$.
\end{lemma}

\begin{proof}
Choose any equivariant embedding $Z \subset T_0(Z)$. Then, if $C =
L$ as a subset of $T_0(Z)$, there is nothing to show. Otherwise
$C$ is a coordinate line of $T_0(Z)$ having a long $T$-weight
$\alpha$, say. If $L \not \subset T_p(Z)$ for a $p \in C^o$ there
is a $T$-eigenfunction $f$ in the ideal of $Z$ in $k[T_0(Z)]$,
such that $df_p(L) \neq 0$. We may assume that $k[T_0(Z)] =
k[z,x_1,x_2,\dots,x_n]$ with $z$ as above corresponding to $L$,
and the $x_i$ corresponding to the long lines of $T_0(Z)$. Then we
write $f = P_0 + P_1z + P_2z^2 +\dots + P_dz^d$ with the $P_i$
$T$-eigenvectors and polynomials in the $x_i$ only. Without loss
of generality $P_iz^i$ has the same weight as $f$. It follows that
$df_p = dP_{0,p} + P_1(p)dz_p$ because $z$ vanishes on $C$. By
assumption $P_1(p)$ is nonzero, implying that there is a monomial
of the form $x^l$ contained in $P_1$, where $x$ is the coordinate
corresponding to $C$ and $l\geq 1$. Thus, the $T$-weight of $f$ is
$l\alpha + \omega$. On the other hand $P_0$ is nonzero. For if
$P_0 = 0$, then $f$ is divisible by $z$, and therefore $f = hz$
for some $h$. But $Z$ is irreducible and clearly $z$ does not
vanish on $Z$, so $h$ vanishes on $Z$. Now $z$ and $h$ vanish in
$p$ forcing $df_p$ to be zero as well, a contradiction.
With $P_0$ being nonzero it follows that there is a monomial in
the $x_i$ of weight $l\alpha + \omega$. This clearly shows that
$\omega = (l\alpha + \omega) - l\alpha$ is contained in the
$\mathbb Z$-submodule of $X(T)$ generated by all long weights of
$T_0(Z)$. The next lemma shows that this is impossible and
therefore ends the proof.

\end{proof}

\begin{lemma}Let $\Gamma$ be a $\mathbb Z$-submodule of $X(T)$
generated by long roots. If the Killing form $F$ is normalized
such that $(\omega,\omega) = 1$ is the short root length, then the
function $f: \Gamma \RA \mathbb Q$ given by $f(\gamma) =
(\gamma,\gamma)$ has actually values in $2\mathbb Z$.
\end{lemma}

\begin{proof}
If $\alpha$, $\beta$ are long roots, then $(\alpha, \beta) \in
\mathbb Z$. Indeed, $(\alpha,\beta) \in \{0,\pm 1, \pm 2\}$ by
general properties of root systems. Hence, $(\gamma,\delta) \in
\mathbb Z$ for all $\gamma, \delta \in \Gamma$, as well. Now
$f(\gamma + \delta) = f(\gamma) + f(\delta) + 2(\alpha,\delta) \in
2\mathbb Z$, if $f(\gamma)$ and $f(\delta)$ are even integers. The
result follows by induction on the length of a shortest
representation $\gamma = \sum n_i \alpha_i$ with $n_i \in \mathbb
Z$ and $\alpha_1,\alpha_2,\dots$ the long generators of $\Gamma$.
The length of such a representation is just $\sum |n_i|$. So, if
$n_1$ is nonzero and positive, then $\gamma = \alpha_1 + (n_1 - 1)
\alpha_1 + \sum_{i > 2} n_i\alpha_i$. The induction hypothesis for
$\alpha_1$ and $(n_1 - 1)\alpha_1 + \sum_{i > 2} n_i \alpha_i$
give the result for $\gamma$ by the above arguments. If $n_1$ is
negative we may use $-\gamma$, since $f(\gamma) = f(-\gamma)$.
Finally, if $n_1$ is zero, we may replace $\alpha_1$ with any
other $\alpha_i$ such that $n_i \not=0$.
\end{proof}

We are now in a position to prove the Theorem.

\begin{proof}[Proof of Theorem \ref{THETA_SUB_TAU}]
We proceed by induction on $\dim Z$ for an irreducible $T$-stable
subvariety $Z \subset X \subset T_x(X)$. Of course there is
nothing to show when $\dim Z \leq 1$. If $\dim Z > 1$, let $L
\subset \Theta_0(Z)$ be any $T$-stable line that has a short
weight $\omega$, say. Let $z$ be a corresponding function of
$k[Z]$. Suppose there is another line with short weight in
$T_0(Z)$. By the previous lemma, if $f$ is a corresponding
function $z$ does not vanish on the tangent cone of $\mathcal
V(f)$. Thanks to Lemma \ref{THETA_SUM}, $z$ does not vanish on the
tangent cone of at least one irreducible component $Z^\prime$ of
$\mathcal V(f)$. In particular this implies that $L$ is contained
in $\Theta_0(Z^\prime)$. By induction $L \subset \tau(Z^\prime,0)
\subset \tau(Z,0)$. This concludes the case that there is a short
line in $T_0(Z)$ different from $L$.
So suppose $L$ is the only line in $T_0(Z)$ with a short weight.
Then $L \subset T_p(Z)$ for all $p \in C^o$ and any curve $C \in
E(Z,0)$. For each such $C$ it then follows that $L \subset
\tau_C(Z,0)$.
By Theorem \ref{TC=TL} all the lines in $\Theta_0(Z)$ with long
$T$-weights are tangent to $T$-curves, so they are contained in
$\tau(Z,0)$.
\end{proof}

We complete this section with an example that shows the $G_2$-restriction is
necessary. We need the following general fact about $\tx$ proved in \cite{car}.
\begin{proposition}  Suppose $X$ is a \SV\ in $G/B$ and $x\in X^T$. Let
$\mathcal H$ 
denote the convex hull in $\Phi \otimes \R$ of the $T$-weights of $TE(X,x)$.
Then every $T$-weight of $\tx$ lies in $\mathcal H$.  
\end{proposition}
\begin{example}\label{G2eg}  Now suppose $\A$ and $\B$ are the short and long simple roots 
in the root system of $G_2$, and consider the \SV\ $X$ in $G_2/B$ corresponding to
$w=r_\B r_\A r_\B r_\A$. By \cite[p. 168]{bl}, the singular locus of $X$ is the \SV\ $X(r_\B r_\A)$,
so $x=r_\B r_\A$ is a maximal singularity. 
By a direct check, the $T$-weights of $TE(X,x)$ are $-\A, \ \B,\ \A+\B$ and $-\LA,$
where $\LA=3\A +2\B$ is the highest root. 
Thus the weights in $\H$ are $-\A,\ \B,\ \A+\B,\ -(\B +2\A)$ and $-\LA,$
The good $T$-curves in $E(X,x)$ correspond to $-\A$ and $-\LA$.
We claim that $-(3\A+\B)$ is a weight in $\tau_C(X,x)$, where $C$ corresponds to
$-\LA$. Indeed, put 
$y=r_{\LA}x$. Then one sees that the weights of $TE(X,y)$ are $\B,\ \A+\B,\
-(\B +2\A)$ and $\LA$. By inspection, $TE(X,y)$ is a $\fg_{-\LA}$-
submodule of $T_y(G_2/B)$, so, by the algorithm
in \cite[\S 8]{ck} (summarized in  Remark \ref{remalg} below), 
the weights of $\tau_C(X,x)$ are obtained by reflecting the
weights of $TE(X,y)$ by $r_{\LA}$. Thus $\tau_C(X,x)$
has weights 
$$r_\LA(\B)=-(3\A +\B),\ r_\LA(\A +\B)= -(2\A +\B), \ r_\LA(-(2\A +\B))=\A +\B,\ 
\mathrm{and} \ r_\LA(\LA)=-\LA.$$
Since $-(3\A +\B)$ isn't in $\H$, Theorem \ref{SUM_TAU} may fail without the $G_2$-restriction.
\end{example}



\section{Proof of Theorems \ref{THB2PAIR} and \ref{BGALG}} \label{PTSV}


The goal of this section is to study the $T$-weights
in $\tau_C(X,x)$ for a \SV \ in $G/B$
and to eventually prove Theorems \ref{THB2PAIR} and \ref{BGALG}. 
As usual, we will suppose throughout that
$G$ does not contain any $G_2$-factors. Let $X=X(w)$  and
assume  $C$ is a good $T$-curve in $X$ \st $C^T=\{x,y \}$,
where $y>x$. Thus we can write $C=\overline{U_{-\mu}x}$, where
$\mu >0$, and it follows that $y=r_\mu x>x$. 
Since $\tau_C(X,x)\subset TE(X,x)$ if $\mu$ is short,
we can ignore this case and suppose $\mu$ is long.
Recall also that if ${\fg}_{\g}\subset\tx$ and $\g$ is long,
then  ${\fg}_{\g}\subset TE(X,x)$.

To begin, we need  a result  similar to Theorem \ref{TC=TL} for $\tau_C(X,x)$.
\begin{lemma}\label{root} Suppose $\g$ is a short root such that 
${\fg}_{\g}\subset
\tau_C(X,x)$. If ${\fg}_{\g}\not \subset TE(X,x)$, then
there exists a long root $\phi$ orthogonal to $\mu$
such that  $\fg_{-\phi}\subset TE(X,x)$, and
\begin{equation}\label{shortroot}\g=-\frac{1}{2}(\mu +\phi).
\end{equation}
In addition, the roots $\g,\mu,\phi$  lie in a copy of $B_2$ contained in 
$\Phi$.
When ${\fg}_{\g}\subset TE(X,x)$, there exists a  $T$-surface
in $S\subset X$  containing $C$ and the $T$-curve 
corresponding to $\g$.
\end{lemma}
\begin{proof} This follows from \cite[Lemma 5.1 and Proposition 5.2]{ck}.

\end{proof} \Black
We will see below that  if ${\fg}_{\g}\not\subset TE(X,x)$, then $\phi >0$.
The notion of an orthogonal
$B_2$-pair arises from the following illuminating example worked out in 
detail in \cite[Example 8.4]{ck}.

\begin{example}\label{B2} Let $G$ be of type $B_2$, and let $w=r_\A r_\B 
r_\A$,
where $\A$ is the short simple root and $\B$ is the long simple
root. Put $X=X(w)$. The singular set  of $X$ is $X(r_\A )$, so $x=r_\A$ 
is $X$'s unique maximal singular point.
There are two good $T$-curves at $x$, namely
$C=\overline{U_{-\B}x}$ and $D=\overline{U_{-(2\A+\B)}x}$.
Suppose $y=r_\B x$ and
$z=r_{2\A+\B} x$. Then
$$T_y(X)=\fg_{-\A} \oplus \fg_{\A+\B}\oplus \fg_\B \quad\quad 
\mathrm{and} \quad\quad
T_z(X)=\fg_{\A} \oplus \fg_{-(\A+\B)}\oplus \fg_{2\A+\B}.$$
Thus, by the algorithm of \cite[\S 3]{ck}, 
$$\tau_C(X,x)=\fg_{\A} \oplus \fg_{-(\A+\B)}\oplus \fg_{-\B}
\quad\quad \mathrm{and} \quad\quad
\tau_D(X,x)=\fg_{\A} \oplus \fg_{-(\A+\B)}\oplus \fg_{-(2\A+\B)}.$$
Note that the weight at $x$ that does not give a $T$-curve,
namely $-(\A+\B)$, is in both Peterson translates. The next result extends
this example to the general case.
\end{example}

\begin{remark}\label{remalg} We will use the algorithm in \cite[\S 3]{ck} in several places
to compute a \pt\ $\tau_C(X,x)$. Let us briefly summarize how this works.
Suppose $C=\overline{U_{-\mu}x}$, where $\mu >0$ and $y=r_\mu x$. 
Consider the weights of the form $\nu +k\mu$ in $T_y(X)$, and form
a (possibly partial) $\mu$ string consisting
of roots of the form $\kappa -j\mu$, where $0\le j\le r,$
\st $y^{-1}(\kappa -j\mu)<0$ for each $j$, but $y^{-1}(\kappa -(r+1)\mu)>0.$
Then the roots $r_\mu (\kappa -j\mu)$ occur as weights in $\tau_C(X,x)$,
and every weight occuring in $\tau_C(X,x)$ arises in this way. 
\end{remark}

Recall that $(\, ,  \, )$ is a $W$-invariant inner product on 
$X(T)\otimes \R$.
Assuming $\g$ is as in the last Lemma, we now say more about $\fg_\g$.
\begin{theorem}\label{phi} 
Suppose  $\g$ is a short root such that ${\fg}_{\g}\subset
\tau_C(X,x)$. If either $(\g,\mu)\ge 0$, or in the
equation $(\ref{shortroot})$ one has $\phi<0$, then  ${\fg}_\g \subset 
TE(X,x)$.
On the other hand, if ${\fg}_\g\not\subset TE(X,x)$, then the following
statements  hold:

\begin{itemize}

\item[(a)] $\g <0$,

\medskip
\item[\rm{(b)}] $(\g,\mu)<0$, hence $\delta:=\g +\mu \in \Phi$,

\medskip
\item[\rm{(c)}] if $x^{-1}(\delta)<0$, then $\fg_\delta \subset 
\tau_C(X,x)\cap TE(X,x)$
$($and, of course, conversely$)$, and

\medskip
\item[\rm{(d)}] $\phi>0$.

\end{itemize}

\end{theorem}

\begin{remark} \label{B2rem}  Example \ref{B2} shows that one can have
${\fg}_{\g}\subset\tau_C(X,x)\cap TE(X,x)$
yet still have $(\g,\mu)<0$.
\end{remark}

\begin{proof}
If $(\g,\mu)\ge 0$, it follows immediately from Lemma \ref{root} that 
$\fg_\g \subset TE(X,x)$.
Suppose $\g$ and  has the form (\ref{shortroot}), where  $\phi<0$,
and put $\D=\g+\mu$. Since $(\g, \mu)<0$, $\delta\in \Phi$.
Moreover, since $\phi <0$, we have $\D>0$. Now if $\g>0$, then $r_\g x<x$,
since $x^{-1}(\g)<0$. Thus $\fg_\g \subset TE(X,x)$ if $\g>0$.

Next, suppose  $\g<0$. We will consider the two cases $x^{-1}(\D)<0$
and $x^{-1}(\D)>0$ separately.
Assume first that $x^{-1}(\D)<0$. Since $\tau_C(X,x)$ is a $\fg 
_\mu$-submodule
of $T_x(X)$ (cf. \cite[\S 3]{ck}) and $\fg_\g \subset \tau_C(X,x)$, we 
therefore know that
$$\fg_\D \oplus \fg_\g \subset \tau_C(X,x).$$
Since $\mu$ is long and there are no $G_2$-factors, Proposition 8.1
\cite{ck} implies
$${\fg}_\D\oplus {\fg}_{\g}\subset T_y(X).$$
Since $\g<0$, we therefore get the inequality $y<r_\g y\leq w$, and 
hence $X$ is also nonsingular
at $r_\g y$.
Moreover, since $\phi <0$ and $x^{-1}(\phi)=y^{-1}(\phi)>0$, it also follows
that $ {\mathbf g}_{-\phi} \subset TE(X,y)$,
which equals $T_y(X)$ since $X$ is smooth  at $y$.
Since there are no $G_2$ factors, $\mu,\ \D,\ -\phi$ constitute
a complete $\g$-string  occuring as $T$-weights of $T_y(X)$.
Letting $E$ be the
good $T$-curve in $X$ such that $E^T=\{y,r_\g y \}$, we have
$\tau_E(X,y)=T_y(X)$, so the string $\mu,\ \D,\ -\phi$ also has to occur
as $T$-weights of $T_{r_\g y}(X)$. In particular, ${\mathbf g}_{-\phi}
\subset  TE(X, r_\g y)=T_{r_\g y}(X)$, and hence $r_\phi r_\g
y\leq  w$. But this means
$$r_\g x=r_\g r_\mu y=r_\g r_\mu r_\g r_\g y=r_\phi r_\g y\leq w,$$
so ${\mathbf g}_\g \subset TE(X,x)$.

Next, assume $x^{-1}(\D)>0$. Since $\mu$ is long, $r_\mu (\D)=\D -\mu=\g$,
hence $y^{-1}(\D)=x^{-1}(\g)<0.$ Thus, since $\D>0$, $\fg_\D \subset 
T_y(X)$.
Furthermore,
$$y^{-1}(-\g)=-x^{-1} r_\mu (\g)=-x^{-1}(\D)<0,$$
so $\fg_{-\g} \subset T_y(X)$.
It follows that $r_\g y <y$. As $-\phi>0$, $U_{-\phi}r_\g y \subset X$
as well. We claim $U_{-\phi}r_\g y \ne r_\g y$, which then proves that
$r_\phi r_\g y \le w$. But
$$(r_\g y)^{-1}(-\phi)=y^{-1}(r_\g(-\phi))=y^{-1}(\mu)<0,$$
hence we get the claim. Finally, we note that $r_\phi r_\g r_\mu =r_\g$,
so it follows that $r_\g x\le w$. Therefore, if $\phi <0$, we get 
$\fg_\g \subset TE(X,x)$.

Now suppose $\fg _\g \not\subset TE(X,x)$.
Then (a) is immediate and (b) follows from (\ref{shortroot}). Since 
$\tau_C(X,x)$
is a $\fg_\mu$-submodule of $T_x(X)$,  $\fg_\delta\subset \tau_C(X,x)$ 
since $x^{-1}(\delta)<0$.
Then $\g$ is given by (\ref{shortroot}), so $(\delta, \mu)\ge 0$
(since $\mu$ is long). Thus, Lemma \ref{root} implies 
$\fg_\delta\subset TE(X,x)$.
On the other hand, if $x^{-1}(\delta)>0$, then $\fg_\delta\not \subset 
T_x(X)$.
This establishes (c). The assumption that $\fg _\g \not\subset TE(X,x)$
immediately implies that $\phi$ is positive giving (d).


\end{proof}


\begin{remark} Let $X$ be a \SV, and suppose $x\in X^T$ is a maximal 
singularity
where $|E(X,x)|=\dim X$.
In this case, the second author has shown that the multiplicity
$\tau_x(X)$ of $X$ at $x$ is exactly $2^d$, where
$$d=|\{ \A\in x(\Phi^-) \mid \fg_\A \subset \tau_C(X,x) 
~~~\mathrm{and}~~ r_\A x \not < w\},$$
for {\em any } good $C\in E(X,x)$ (\cite{kut}).
\end{remark}

\begin{theorem} \label{B2PAIRTH} 
Suppose $C=\overline{U_{-\mu}x}$ is a good $T$-curve, where
$\mu >0$, and let $y=r_\mu x$. 
Assume
$\fg _\g \subset \tau_C(X,x)$ but $\fg _\g \not\subset \TX$. Then there
exists a positive root $\phi$ \st $\{\mu,\phi\}$ is an orthogonal 
$B_2$-pair for $X$ at $x$
\st $\g=-1/2(\mu +\phi)$.
Conversely, suppose that for some $\phi >0$, $\{\mu,\phi\}$ is an orthogonal 
$B_2$-pair for $X$ at $x$, and $\g=-1/2(\mu +\phi)$. Then $\fg_\g \subset \tau_C(X,x)$.
\end{theorem}

\begin{proof} Suppose $\fg _\g \subset \tau_C(X,x)$ but $\fg _\g 
\not\subset \TX$.
By Lemma \ref{root} and Theorem \ref{phi}, there exists a long positive 
root $\phi$ orthogonal to $\mu$
\st $\g=-1/2(\mu +\phi)$. Put $y=r_\mu x$, and note $X$ is smooth at $y$.
To show that $\{\mu,\phi\}$ is an orthogonal  $B_2$-pair,  we have to 
consider two cases.

\ms {\bf Case 1}. $\mu$ is simple. Then $\A=\g+\phi$ is the short simple 
root.
We have to show that if $\fg_\g \not\subset \TX$, then $r_\A x<x$ and
$r_\A r_\mu x\le w.$ But $\fg_\g \not\subset \TX$ implies
$x^{-1}(\A)<0$, since if $x^{-1}(\A)>0$, then the fact that
$\g=-\phi +\A$ would say  $\fg_\g \subset \TX.$ Hence $r_\A x<x$.

Since $r_\mu (\A)=-\g$, it follows that $y^{-1}(\A)=x^{-1}(-\g)>0,$
so $\fg_{-\A}\subset T_y(G/B)$. But $y^{-1}(\g)=x^{-1}(-\A)>0$,
hence $\fg_\g \not\subset T_y(G/B)$. Hence,
by the algorithm for computing
the \pt\ in \cite{ck} and the fact that $\fg _\g \subset \tau_C(X,x)$,
we infer that $\fg_{-\A} \subset T_y(X)$. Therefore, $r_\A y=r_\A r_\mu 
x\le w$,
as was to be shown.

\ms {\bf Case 2}. $\phi$ is simple. Here $\A=\g+\mu$ is the short simple 
root,
and $r_\mu (\g)=\A$. As in Case 1, $x^{-1}(\A)<0$, so $r_\A x<x$.
Now $y^{-1}(\A)=x^{-1}(\g)<0,$ so $r_\A y<y$ and hence $\fg_\A \subset 
T_y(X)$.
Also, $y^{-1}(\g)=x^{-1}(\A)<0$, so $\fg_\g \subset T_y(G/B)$.
Thus the algorithm for $\tau_C(X,x)$ says that $\fg_\A \subset \tau_C(X,x)$.
But as $\fg_\g \subset \tau_C(X,x)$ too, we have to conclude that
$\fg_\g \subset T_y(X)$, due to the fact that $\g$ and $\A$ comprise a 
$\mu$ string.
Hence $r_\g y \le w$. But since we are in a $B_2$ where $\A$ and $\phi$ 
are the simple
roots, $r_\g r_\mu =r_\A r_\phi$. Hence $r_\A r_\phi x \le w$, so  Case 2 is finished.

To prove the converse, we need to consider Cases 1 and 2 again
with the assumption that $x^{-1}(\A)<0$, which follows from the 
condition that $r_\A x <x$. The argument is, in fact, very similar to 
the above, but
we will outline it anyway. Assume first that $\mu=\B$, i.e. $\mu$ is simple.
As $r_\A r_\B x \le w$, we see that $r_\A y \le w$. But $y^{-1}(-\A)=x^{-1}(\g)<0$,
consequently $\fg_{-\A} \subset T_y(X)$. Also, $y^{-1}(\g)=
x^{-1}(-\A)>0$, so
$\fg_\g \not\subset T_y(G/B)$. Thus by the algorithm for computing
$\tau_C(X,x)$, the weight $r_{\B}(-\A)$ occurs in $\tau_C(X,x)$. 
Hence, $\fg_\g \subset  \tau_C(X,x)$.

On the other hand, if $\phi$ is simple, then $\mu =\B +2\A$.
Thus,  $y^{-1}(\A)=x^{-1}(\g)<0$, so $r_\A y <y$, hence $\fg_\A \subset T_y(X)$.
But $r_\A r_\phi x \le w$ means $r_\g r_\mu x \le w$,
that is, $r_\g y \le w$. As $y^{-1}(\g)=x^{-1}(\A)<0$,
$\fg_\A +\fg_\g \subset T_y(Y)$. Since $\A$ and $\g$ make up
a $\B +2\A$-string in $B_2$,  $\fg_\A +\fg_\g \subset \tau_C(X,x)$ also. This
finishes the proof.
\end{proof}

We now prove Theorems \ref{THB2PAIR} and \ref{BGALG}. 
\begin{proof}[Proof of Theorem \ref{THB2PAIR}] 
Suppose $\fg_\g \subset \tx$. Since $x$ is either smooth
or a maximal singularity,  Theorem \ref{SUM_TAU}
$\fg_\g \subset \tau_C(X,x)$ for some good $C$. If $C$ is short, then
$\tau_C(X,x)\subset TE(X,x)$, by Theorem \ref{TH1}, hence $\tau_C(X,x)\subset \TX$.
Thus we can suppose $C$ is long. But then, by Theorem \ref{B2PAIRTH},
either $\fg_\g \subset \TX$ or there exists a $B_2$-pair $\{\mu,\phi\}$ for $X$ at $x$
\st $\g=-1/2(\mu +\phi)$. Hence Theorem \ref{THB2PAIR} is proven.
\end{proof}

\begin{proof}[Proof of Theorem \ref{BGALG}] Suppose $C\in E(X,x)$ 
is good and $\dim TE(X,x)=\dim \TX=\dim X$. If $C$ is short, then $X$ is smooth at $x$
by Theorem \ref{TH1}. Hence we may suppose $C$ is long.
Suppose there  exists a $T$-line $\fg_\g$ in $\tau_C(X,x)$ which is not 
in $\TX$.
Then by Theorem \ref{B2PAIRTH}, there is an orthogonal $B_2$-pair 
$\{\mu,\phi\}$
for $X$ at $x$ for which $\g=-1/2 (\mu +\phi)$. But then by assumption,
$\fg_\g \subset TE(X,x)$.
This contradicts the choice of $\fg_\g$, so $\tau_C(X,x)\subset \TX=TE(X,x)$. Hence,
by Theorem \ref{TH1} again, $X$ is smooth at $x$.

For the converse, suppose $X$ is smooth at $x$. Then conditions (1) and (2) of Theorem
\ref{THB2PAIR} clearly hold. Suppose $\{\mu,\phi\}$ is a $B_2$-pair
for $X$ at $x$ and  $\g=-1/2(\mu +\phi)$. By the converse assertion
of Theorem \ref{B2PAIRTH}, $\fg_\g \subset \tau_C(X,x)$, where $C\in E(X,x)$ is the
$T$-curve of weight $\mu$ at $x$. Since $x$ is smooth, $\tau_C(X,x)=TE(X,x)$,
so $\fg_\g \subset TE(X,x)$. 
\end{proof}


\newpage

{\footnotesize
\begin{theckbibliography}{00}

\bibitem{bb} S.\ Billey and T.\ Braden: {\em Lower Bounds for Kahzdan-Lusztig Polynomials from Patterns},
To appear in Transf. Groups.

\bibitem{bl} S.\ Billey and V.\ Lakshmibai: {\it Singular loci of Schubert
varieties}. Progress in Mathematics {\bf 182}, Birkh\"auser Boston-
Basel-Berlin, 2000.

\bibitem{bp}  S.\ Billey and A.\ Postnikov: {\em Smoothness
of Schubert Varieties Via Patterns in 
Root Systems} (2003), arXiv:math.CO/0205179 v1.


\bibitem{bg} B.\ Boe and W.\ Graham: {\em A lookup conjecture for
rational smoothness}. Amer. J. Math.

\bibitem{cp}  J.\ Carrell: {\em The Bruhat Graph of a Coxeter Group, a
Conjecture of Deodhar, and Rational Smoothness of Schubert Varieties}.
Proc. Symp. in Pure Math. A.M.S. {\bf 56} (1994), Part I, 53-61.

\bibitem{car} J.\ Carrell: {\em The span of the tangent cone of a
  Schubert variety}. Algebraic Groups and Lie Groups, Australian Math.
Soc. Lecture Series {\bf 9}, Cambridge Univ. Press (1997), 51-60.


\bibitem{ck}  J.\ Carrell and J. \ Kuttler: {\em Singular points of 
$T$-varieties in $G/P$
and the Peterson map}. Invent. Math. {\bf 151}  (2003), 353-379
DOI:10.1007/s00222-002-0256-5.

\bibitem{chev} C.\ Chevalley: {\it Sur les d\' ecompositions cellulaires des
espaces $G/B$}, Proc. Symp.2 Pure Math. A.M.S. {\bf 56} (1994), Part I, 
1-25.


\bibitem{kl} D.\ Kazhdan and G.\ Lusztig:  {\em Representations of 
Coxeter groups and Hecke algebras}, Invent. Math. {\bf 53} (1979),
165-184.

\bibitem{ku}  S.\ Kumar: {\it Nil Hecke ring and singularity of Schubert
varieties}, Inventiones Math. {\bf 123} (1996), 471-506.

\bibitem{kut} J.\ Kuttler: {\it The Singular Loci of $T$-Stable
Varieties in $G/P$}, thesis, U. of Basel (2003)









\end{theckbibliography} }

\bigskip
\noindent
{\tiny
\begin{tabular}{l}
James B.\ Carrell\\
Department of Mathematics\\
University of British Columbia\\
Vancouver, Canada V6T 1Z2\\
carrell$@$math.ubc.c{a}
\end{tabular}}

\bigskip
\noindent
{\tiny
\begin{tabular}{l}
Jochen Kuttler\\
Mathematisches Institut\\
Universit\"at Basel\\
CH-4051 Basel\\
Switzerland\\
kuttler$@$math.unibas.c{h}
\end{tabular}}

\end{document}